\magnification=\magstep1
\hsize=16.5 true cm 
\vsize=25 true cm
\font\bff=cmbx10 scaled \magstep1

\font\bffg=cmbx10 scaled \magstep3

\font\smc=cmcsc10 
\parindent0cm
\overfullrule=0cm
\def\cl{\centerline}           %
\def\rl{\rightline}            %
\def\bp{\bigskip}              %
\def\mp{\medskip}              %
\def\sp{\smallskip}            %
           %
\def\Bbb#1{\hbox{\boldmas #1}} %
\centerline{\bffg On the variety of Euclidean point sets} 
\bigskip\sp
\centerline{\bff Gerald Kuba}
\bigskip\bp

{\expandafter\edef\csname amssym.def\endcsname{%
       \catcode`\noexpand\@=\the\catcode`\@\space}

\catcode`\@=11

\def\undefine#1{\let#1\undefined}
\def\newsymbol#1#2#3#4#5{\let\next@\relax
 \ifnum#2=\@ne\let\next@\msafam@\else
 \ifnum#2=\tw@\let\next@\msbfam@\fi\fi
 \mathchardef#1="#3\next@#4#5}
\def\mathhexbox@#1#2#3{\relax
 \ifmmode\mathpalette{}{\m@th\mathchar"#1#2#3}%
 \else\leavevmode\hbox{$\m@th\mathchar"#1#2#3$}\fi}
\def\hexnumber@#1{\ifcase#1 0\or 1\or 2\or 3\or 4\or 5\or 6\or 7\or 8\or
 9\or A\or B\or C\or D\or E\or F\fi}

\font\tenmsa=msam10
\font\sevenmsa=msam7
\font\fivemsa=msam5
\newfam\msafam
\textfont\msafam=\tenmsa
\scriptfont\msafam=\sevenmsa
\scriptscriptfont\msafam=\fivemsa
\edef\msafam@{\hexnumber@\msafam}
\mathchardef\dabar@"0\msafam@39
\def\dashrightarrow{\mathrel{\dabar@\dabar@\mathchar"0\msafam@4B}}
\def\dashleftarrow{\mathrel{\mathchar"0\msafam@4C\dabar@\dabar@}}

\def\ulcorner{\delimiter"4\msafam@70\msafam@70 }
\def\urcorner{\delimiter"5\msafam@71\msafam@71 }
\def\llcorner{\delimiter"4\msafam@78\msafam@78 }
\def\lrcorner{\delimiter"5\msafam@79\msafam@79 }
\def\yen{{\mathhexbox@\msafam@55}}
\def\checkmark{{\mathhexbox@\msafam@58}}
\def\circledR{{\mathhexbox@\msafam@72}}
\def\maltese{{\mathhexbox@\msafam@7A}}

\font\tenmsb=msbm10
\font\sevenmsb=msbm7
\font\fivemsb=msbm5
\newfam\msbfam
\textfont\msbfam=\tenmsb
\scriptfont\msbfam=\sevenmsb
\scriptscriptfont\msbfam=\fivemsb
\edef\msbfam@{\hexnumber@\msbfam}
\def\Bbb#1{{\fam\msbfam\relax#1}}
\def\widehat#1{\setbox\z@\hbox{$\m@th#1$}%
 \ifdim\wd\z@>\tw@ em\mathaccent"0\msbfam@5B{#1}%
 \else\mathaccent"0362{#1}\fi}

\def\widetilde#1{\setbox\z@\hbox{$\m@th#1$}%
 \ifdim\wd\z@>\tw@ em\mathaccent"0\msbfam@5D{#1}%
 \else\mathaccent"0365{#1}\fi}
\font\teneufm=eufm10
\font\seveneufm=eufm7
\font\fiveeufm=eufm5
\newfam\eufmfam
\textfont\eufmfam=\teneufm
\scriptfont\eufmfam=\seveneufm
\scriptscriptfont\eufmfam=\fiveeufm

\newsymbol\risingdotseq 133A
\newsymbol\fallingdotseq 133B
\newsymbol\complement 107B
\newsymbol\nmid 232D
\newsymbol\rtimes 226F
\newsymbol\thicksim 2373

\font\eightmsb=msbm8   \font\sixmsb=msbm6   \font\fivemsb=msbm5
\font\eighteufm=eufm8  \font\sixeufm=eufm6  \font\fiveeufm=eufm5
\font\eightrm=cmr8     \font\sixrm=cmr6     \font\fiverm=cmr5
\font\eightbf=cmbx8    \font\sixbf=cmbx6    
      \font\eighti=cmmi8   \font\sixi=cmmi6
\font\ninesy=cmsy9     \font\eightsy=cmsy8  \font\sixsy=cmsy6
     \font\eightit=cmti8  
     \font\eightsl=cmsl8  
     \font\eighttt=cmtt8

\font\eightsmc=cmcsc8
\newskip\ttglue
\newfam\smcfam
\def\eightpoint{\def\rm{\fam0\eightrm}%
  \textfont0=\eightrm \scriptfont0=\sixrm \scriptscriptfont0=\fiverm
  \textfont1=\eighti \scriptfont1=\sixi \scriptscriptfont1=\fivei
  \textfont2=\eightsy \scriptfont2=\sixsy \scriptscriptfont2=\fivesy
  \textfont3=\tenex \scriptfont3=\tenex \scriptscriptfont3=\tenex
  \def\smc{\fam\smcfam\eightsmc}
  \textfont\smcfam=\eightsmc          
\textfont\eufmfam=\eighteufm              \scriptfont\eufmfam=\sixeufm
     \scriptscriptfont\eufmfam=\fiveeufm
\textfont\msbfam=\eightmsb            \scriptfont\msbfam=\sixmsb
     \scriptscriptfont\msbfam=\fivemsb
\def\it{\fam\itfam\eightit}%
  \textfont\itfam=\eightit
  \def\sl{\fam\slfam\eightsl}%
  \textfont\slfam=\eightsl
  \def\bf{\fam\bffam\eightbf}%
  \textfont\bffam=\eightbf \scriptfont\bffam=\sixbf
   \scriptscriptfont\bffam=\fivebf
  \def\tt{\fam\ttfam\eighttt}%
  \textfont\ttfam=\eighttt
  \tt \ttglue=.5em plus.25em minus.15em
  \normalbaselineskip=9pt
  \def\MF{{\manual opqr}\-{\manual stuq}}%
  \let\big=\eightbig
  \setbox\strutbox=\hbox{\vrule height7pt depth2pt width\z@}%
  \normalbaselines\rm}
\def\eightbig#1{{\hbox{$\textfont0=\ninerm\textfont2=\ninesy
  \left#1\vbox to6.5pt{}\right.\n@space$}}}

\catcode`@=13 

\vbox{\eightpoint
{\bf Abstract.} We construct  
a continuum of non-homeomorphic compact subspaces of ${\Bbb R}$
without singleton components. Thus 
from the purely topological point of view 
the real line ${\Bbb R}$
contains not only more closed sets than open sets but also more
closures of open sets than open sets. 
On the other hand, we show that this discrepancy vanishes      
either if the topological point of view is sharpened in 
the metrical or in the order-theoretical direction, or if 
${\Bbb R}$ is replaced with ${\Bbb R}^n$ for $\,n\geq 2\,$.
Furthermore, we track down a continuum
of topological types of closed and totally disconnected subsets of ${\Bbb R}$.
In doing so we also track down a continuum of metrical types 
of infinite, discrete subsets of $[0,1]$. 
(As a consequence, any countably infinite discrete space has
a continuum of non-homeomorphic metrizable compactifications.)}} 

\bp\mp
{\bff 1. Unions of intervals}
\mp
 Let $\,{\cal U}\,$ be the family of all
open subsets of the Euclidean real line $\,{\Bbb R}\,$.
Motivated by the fact that each open subset of $\,{\Bbb R}\,$ is a union of 
mutually disjoint open intervals, we define a family 
$\,{\cal A}\,$ of point sets so that $\;X\in{\cal A}\;$ if and only if
$\,X\,$ is a closed subset of $\,{\Bbb R}\,$ which 
is a union of mutually disjoint nondegenerate closed intervals.
In other words, the members of $\,{\cal A}\,$ are precisely the closed
subspaces of $\,{\Bbb R}\,$ where no component of any space is 
a singleton.
\medskip
It is common opinion that the 
topological structure of an arbitrary closed subset 
$\,{\Bbb R}\,$ may be more 
complicated than of any open subset of $\,{\Bbb R}\,$,
although grounds for this opinion are rather informal.
Our first goal is to support this view by 
pointing out that the structural discrepancy in question is revealed 
by a clear cardinal discrepancy already when only 
point sets in the family $\;{\cal U}\cup{\cal A}\;$
are considered. In fact, there are more topological types
of members of  $\,{\cal A}\,$ than of members of $\,{\cal U}\,$!
\smallskip
Naturally, both families 
$\,{\cal U}\,$ and $\,{\cal A}\,$ have the cardinality 
$\,c\,$ of the continuum ${\Bbb R}$. But the family $\,{\cal U}\,$ 
contains only countably many topologically distinct members. 
Indeed, since each $\,U\in{\cal U}\,$ can be written as 
a union of countably many 
mutually disjoint open intervals, if 
$\;\emptyset\not= U\in{\cal U}\;$ then 
there is precisely one $\;n\in{\Bbb N}\cup\{\infty\}\;$
so that $\,U\,$ and 
$\;\bigcup\limits_{k=1}^n ]k,k+1[\;$
are homeomorphic subspaces of $\,{\Bbb R}\,$.
(In order to avoid potential misinterpretations,
$\,0\not\in{\Bbb N}\,$, i.e.~$\;{\Bbb N}\,=\,\{1,2,3,...\}\,$.)
In particular, each open subspace of $\,{\Bbb R}\,$ with infinitely many 
components is homeomorphic to 
$\,{\Bbb R}\setminus{\Bbb Z}\,$.
On the other hand, the following theorem shows that 
there are $\,c\,$ topologically distinct point sets in
the family $\,{\cal A}\,$. 
\medskip
{\bf Theorem 1.} {\it There are $\,c\,$ mutually non-homeomorphic 
compact subspaces of $\,{\Bbb R}\,$
without singleton components.}
\medskip
The situation is different when, instead of topological types
of point sets in $\,{\Bbb R}\,$, 
{\it metrical} types are considered, which means that
{\it continuity} is sharpened to {\it uniform continuity}. 
(The metric is the inherited Euclidean metric of $\,{\Bbb R}\,$.)
A fortiori, topologically distinct point sets are always
metrically distinct. Thus the interiors of the $\,c\,$
topologically distinct compact point sets given by Theorem 1
must be metrically distinct because every $\,A\in{\cal A}\,$ 
equipped with the Euclidean metric 
is a completion of the interior of $\,A\,$ 
equipped with the Euclidean metric.  
Therefore, the 
total number of metrical types of the open point sets in $\,{\Bbb R}\,$
is $\,c\,$ and hence greater than the total number of their 
topological types. (As a consequence,
there exists a collection of $\,c\,$ metrically distinct and 
topologically similar open subsets of $\,{\Bbb R}\,$.)
Certainly, metrically distinct {\it compact} 
subspaces of $\,{\Bbb R}\,$ cannot be homeomorphic.
But, as the following proposition shows
in an illustrative way, 
it is possible to track down $\,c\,$ members of $\,{\cal A}\,$ 
which are metrically distinct and 
topologically similar. And there is also an
illustrative stack of $\,c\,$ 
metrically distinct and 
topologically similar open subsets
of $\,{\Bbb R}\,$.
\eject
\medskip\smallskip
{\bf Proposition 1.} {\it For each real number $\,u\geq 2\,$ define 
$\,X_u\in{\cal A}\,$ and $\,X_u^\circ\in{\cal U}\,$ via
\smallskip
\centerline{$X_u\;:=\;\bigcup\limits_{n=1}^\infty
\big[2^{u^n},2^{u^n}+u^n\big]\quad$ and \quad
$X_u^\circ\;=\;\bigcup\limits_{n=1}^\infty
\big]2^{u^n},2^{u^n}+u^n\big[\;.$}
\smallskip
If $\;2\leq v<w\,$ then there is no uniformly continuous bijection 
from $\,X_v\,$ onto $\,X_w\,$ or from $\,X_v^\circ\,$ onto $\,X_w^\circ\,$.}
\bigskip\medskip
Beside the topological and the metrical view 
there is a third natural way to look at 
the point sets in the families $\,{\cal U}\,$ and $\,{\cal A}\,$.
Two sets $\,X,Y\subset{\Bbb R}\,$ are {\it order-isomorphic}
if and only if there exists a strictly increasing function from 
$\,X\,$ onto $\,Y\,$. Of course, order-isomorphic sets
$\,X,Y\subset{\Bbb R}\,$ need not be homeomorphic subspaces of 
$\,{\Bbb R}\,$.
(Consider for example $\,X\,=\,[2,3]\,$ and 
$\,Y\,=\,\{1\}\cup\,]2,3]\,$.) However, if $\,X,Y\subset{\Bbb R}\,$ are 
open or closed then
the spaces $\,X,Y\,$ must be homeomorphic if the sets 
$\,X,Y\,$ are order-isomorphic. 
Because it is plain that 
the Euclidean topology restricted to a closed or open 
set $\,S\subset{\Bbb R}\,$ coincides 
with the order topology on $\,S\,$ induced by the natural ordering 
of the real numbers in $\,S\,$. And, naturally, any 
order isomorphism between two linearly ordered spaces
is a homeomorphism with respect to their order topologies.
In particular, 
topologically distinct sets in the family $\,{\cal U}\cup{\cal A}\,$
are never order-isomorphic. On the other hand,
the $\,c\,$ metrically distinct open resp.~closed sets 
in Proposition 1 are obviously order-isomorphic. 
It is also possible to  establish
a completely converse situation.
\mp
{\bf Proposition 2.} {\it There are $\,c\,$ metrically (and hence 
topologically) similar sets in the family $\,{\cal U}\,$
and in the family $\,{\cal A}\,$, respectively,  which 
are mutually not order-isomorphic.} 
\mp
Thus, other than concerning topological types
and similarly as concerning metrical types, 
there is no discrepancy between the total numbers  
of the order types of open sets and 
of the order types of closed sets   
in the real line.
\bp\mp
{\bff 2. Unions of cubes}
\mp
The  cardinal discrepancy between all topological types of open 
and all topological 
types of closed sets 
in the realm of linear point sets 
already vanishes in the realm of planar point sets.
Indeed, the following theorem shows that for arbitrary dimensions 
$\,n\geq 2\,$ the Euclidean space
$\,{\Bbb R}^n\,$ contains $\,c\,$ topologically distinct open point sets 
whose closures have no singleton components 
and are topologically distinct as well.
(As usual, $\,\overline{X}\,$ denotes the closure 
of $\,X\,$.)
\medskip\smallskip
{\bf Theorem 2.} {\it For each $\,n\geq 2\,$ 
there is a family $\,{\cal F}_n\,$ of 
open subsets of the Euclidean space $\,{\Bbb R}^n\,$ such that $\,{\cal F}_n\,$ has 
cardinality $\,c\,$ and neither $\,U,V\,$ nor 
$\,\overline U,\overline V\,$ are homeomorphic 
subspaces 
of ${\Bbb R}^n$ whenever $\,U,V\in{\cal F}_n\,$ and $U\!\not=\!V$. 
Moreover, the family ${\cal F}_n$ can be chosen so that 
for every $\,X\in{\cal F}_n\,$ the set $\,\overline X\,$    
is a compact union of closed cubes of the form 
$\,[a_1,a_1+h]\times\cdots\times[a_n,a_n+h]\,$
with $\,h>0\,$ where the interiors of distinct cubes are always disjoint.
Alternatively, the family $\,{\cal F}_n\,$ can be chosen so that 
$\,\overline X\,$ is a union of unit cubes 
$\,[k_1,k_1\!+\!1]\times\cdots\times[k_n,k_n\!+\!1]\,$ with
$\;k_1,...,k_n\in{\Bbb Z}\;$
for every $\,X\in{\cal F}_n\,$.}
\medskip\smallskip
It is impossible that 
every set $\,U\,$ in the uncountable family $\,{\cal F}_n\,$
is a union of mutually disjoint open cubes
(or that $\,\overline{U}\,$ is a union of mutually disjoint 
compact cubes for every $\,U\in{\cal F}_n\,$.) 
Because two open subspaces of $\,{\Bbb R}^n\,$ 
where each component is an open cube are homeomorphic
if and only if the total numbers of components coincide.
But if one sets a value on {\it disjoint} cubes,
it is possible to achieve the following results.
\vfill\eject
\medskip\smallskip
{\bf Proposition 3.} {\it For each of the $\,c\,$ sets $\,S\subset{\Bbb N}\,$ 
define an open set $\,Y_S\subset\,]\!-1,1[^n\,$ via 
\smallskip
\centerline{$Y_S\;:=\;\,]\!-\!1,0[^n\,\cup
\bigcup\limits_{m=1}^\infty \Big(\,\big]2^{-2m},2^{-2m+1}\big[\,
\setminus
\bigcup\limits_{s\in S} 
\big\{\,2^{-2s}\!+\!{k\over s+1}2^{-2s}\;\big|\;k\,=\,1,2,...,s\,\big\}\Big)^n
\;.$}
\smallskip            
Obviously, all $\,Y_S\,$ are unions of infinitely many mutually disjoint 
open cubes
and hence homeomorphic spaces, 
and $\;\overline{Y_S}\,=\,
[-1,0]^n\cup\bigcup_{m=1}^\infty \big[2^{-2m},2^{-2m+1}]^n\;$
for every $\,S\subset{\Bbb N}\,$.
But whenever $\;S\not=S'\,$, there is no bijection $\,f\,$
from $\,Y_S\,$ onto $\,Y_{S'}\,$ such that 
both $\,f\,$ and $\,f^{-1}\,$ are uniformly continuous.}

\medskip\smallskip
{\bf Theorem 3.} {\it For each dimension $\,n\,$ there exists a 
family $\,{\cal V}_n\,$ of open subsets of $\,{\Bbb R}^n\,$
such that {\rm (i)} $\,{\cal V}_n\,$ has the cardinality $\,c\,$; 
{\rm (ii)} each $\,V\in{\cal V}_n\,$ is a union of mutually disjoint open cubes;
{\rm (iii)} $\,\overline{V}\,$ is a union of mutually disjoint 
nondegenerate compact cubes if $\,V\in{\cal V}_n\,$; {\rm (iv)}
all $\,V\in{\cal V}_n\,$ are metrically distinct but topologically similar;
{\rm (v)} all $\,\overline{V}\,(V\in{\cal V}_n)\,$ are
mutually non-homeomorphic compact subspaces of $\,{\Bbb R}^n\,$.}

\bigskip\mp
{\bff 3. Proof of Theorem 1}
\mp
In the following we need {\it Cantor derivatives}
but in order to keep the story simple we use only {\it finite} derivatives.
If $\,P\,$ is a point set in a Hausdorff space then 
the first derivative $\,P'\,$ of 
$\,P\,$ is the set of all limit points of $\,P\,$.
The first derivative of any set is closed. And $\,P\,$ is closed 
if and only if $\,P'\subset P\,$. 
Further, with $\;P^{(0)}=P\,$, for every $\;k\,=\,1,2,3,...\;$
the $k$-th derivative $\,P^{(k)}\,$ of $\,P\,$
is given by $\;P^{(k)}\,=\,(P^{k-1)})'\,.$
Consequently, all derivatives of a closed set $\,A\,$ are closed
and $\;A=A^{(0)}\supset A^{(1)}\supset A^{(2)}\supset A^{(3)}\supset\cdots\,$.
(And possibly but not necessarily, $\,A^{(k)}=A^{(m)}\,$ 
whenever $\,k\geq m\,$
for some $\,m\in{\Bbb N}\,$.)
For abbreviation  
let $\;h(x)\,:=\,{2\over \pi}\arctan x\,$. 
(Then $\,h\,$ is a strictly increasing function from 
$\,[0,\infty[\;$ onto $\,[0,1[\;$.)
\smallskip
In order to prove Theorem 1
we construct a compact subspace $\,X_S\,$ of $\,{\Bbb R}\,$ 
without point components for each 
infinite $\;S\subset{\Bbb N}\;$
so that $\;X_S,X_{S'}\;$ are never homeomorphic for distinct sets $\,S,S'\,$.
Define for each $\;n\in{\Bbb N}\;$ a countable subset 
$\,K_n\,$ of the interval $\,[5n,5n+1]\,$ in the following way.
For arbitrary $\;X\subset [0,1]\;$ define $\;F(X)\subset[0,1]\;$ by
\medskip
\centerline{$\;F(X)\,:=\,
h(\{\,n-1+x\;\,|\,\;n\in{\Bbb N}\;\land\;x\in X\,\})\cup\{1\}\;$}
\medskip
and starting with $\;A_1\,=\,\{\,1-{1\over m}\;|\;m\in{\Bbb N}\,\}\cup\{1\}\;$
put $\;A_{n+1}\,=\,F(A_n)\;$ for $\;n\,=\,1,2,3,...\;$
and define $\;K_n\,:=\,\{\,x+5n\;|\;x\in A_n\,\}\;$
for each $\;n\in{\Bbb N}\,$.
Obviously, $\,K_n\,$ is a closed subset of $\,[5n,5n+1]\,$
with $\;\max K_n\,=\,5n+1\;$ for each $\,n\in{\Bbb N}\,$.
Hence $\,K_n\,$ is always compact.
Furthermore it is evident that $\,K_n\,$ is well-ordered by 
the natural ordering $\,\leq\,$. 
(Besides, one may realize that the order type of $\,(K_n,\leq)\,$ 
is $\,\omega^n+1\,$.)
\medskip
By construction, for each $\,n\in{\Bbb N}\,$ 
the $k$-th Cantor derivative $\,K_n^{(k)}\,$ is infinite 
whenever $\,k<n\,$ and empty whenever $\,k>n\,$ and
$\;K_n^{(n)}=\{5n+1\}\,$.
For each $\;n\in{\Bbb N}\;$ let $\,g_n\,$ be
the reflection in the point $\,5n+2\,$, whence
$\;g_n(x)\,=\,10n+4-x\;$ for $\;x\in{\Bbb R}\;$ and 
$\;g_n([5n,5n+1])=[5n+3,5n+4]\,$.
For every $\;a\in K_n\;$ choose $\;0<\epsilon(a)\leq 1\;$
such that $\;[a,a+\epsilon(a)]\cap K_n\,=\,\{a\}\;$
whenever $\,a\in K_n\,$ 
and put $\;\epsilon(5n+1)=1\,$. 
(For example put $\;\epsilon(a)\,=\,(a'-a)/2\;$ where
$\;a'\,=\,\min\,\{\,x\in K_n\;|\;x>a\,\}\;$
whenever $\;a\,\in\,K_n\setminus\{5n+1\}\,$.) 
Finally, for each infinite set $\;S\subset{\Bbb N}\;$ define 
\medskip
\centerline{$\;X_S\;:=\;h\big(
\bigcup\,\big\{\,\bigcup\,\{\,[a,a+\epsilon(a)]\cup g_n([a,a+\epsilon(a)])\;|
\;a\in K_n\,\}\;\,\big|\,\;n\in S\,\big\}\big)\,\cup\,[1,2]\,$.}
\medskip
It is plain that $\,X_S\,$ is always a closed and hence compact
subset of $\,[0,2]\,$. Obviously, 
all components of the space $\,X_S\,$ are compact (and nondegenerate) intervals
and $\;C_n\,:=\,h([5n+1,5n+3])\;$ is a component of $\,X_S\,$ 
for every $\;n\in S\,$.
Hence we can write $\;X_S\,=\,\bigcup\,\{\,[a_j,b_j]\;|\;j\in{\Bbb N}\,\}\;$ 
where always 
$\;a_j<b_j\;$ and the intervals $\,[a_j,b_j]\,$ 
are mutually disjoint. 
\smallskip
Consider the point set 
$\;B_S\,:=\,\{\,a_j\;|\;j\in{\Bbb N}\,\}\cup\{\,b_j\;|\;j\in{\Bbb N}\,\}\,$,
which clearly is the boundary 
of the point set $\,X_S\,$ in the Euclidean space $\,{\Bbb R}\,$.
The point set $\,B_S\,$ is also topologically determined 
within the subspace $\,X_S\,$ of $\,{\Bbb R}\,$ because
$\,B_S\,$ equals the set of all points $\,x\in X_S\,$
so that for every component $\,C\,$ of the space $\,X_S\,$
the point set $\,C\setminus\{x\}\,$ remains connected
in the space $\,X_S\,$.
Let $\,{\cal L}_S\,$ be the family of all
components $\,C\,$ of $\,X_S\,$ such that $\,C\,$
contains precisely two limit points of $\,B_S\,$.
(Thus the family $\,{\cal L}_S\,$ is also topologically determined 
with respect to the space $\,X_S\,$.)
By construction we have $\;{\cal L}_S\,=\,\{\,C_n\;|\;n\in S\,\}\;$
for each infinite $\,S\subset{\Bbb N}\,$.
(Any component $\,C\,$ of $\,X_S\,$ with $\,C\not=C_n\,$ 
for every $\,n\in S\,$ contains at most one limit point of $\,B_S\,$.)
Moreover, if $\,S\,$ is any infinite subset of $\,{\Bbb N}\,$ 
and if $\,k\in{\Bbb N}\,$ and $\;n\in S\;$
then $\;B_S^{(k)}\cap C_n\,=\,\emptyset\;$
when $\;k>n\;$ and $\;B_S^{(k)}\cap C_n\,\not=\,\emptyset\;$
when $\;k\leq n\,$. Consequently, 
\medskip
\centerline{$S\;=\;\big\{\,\min\,\{\,m\in{\Bbb N}\;|\;B_S^{(m+1)}\cap C=\emptyset\,\}
\;\,\big|\;\,C\in{\cal L}_S\,\big\}$}
\medskip
and hence the set $\,S\,$ is completely determined 
by the topology of the space $\,X_S\,$.
Thus for distinct infinite sets $\;S,S'\subset{\Bbb N}\;$
the spaces $\;X_S,X_{S'}\;$ cannot be homeomorphic.
\bigskip
{\it Remark.} The clue in the previous proof is to approximate
certain intervals from {\it both} the left and the right.
The proof would not work with approximations, say, from the left.
Because if we consider the compact spaces
\smallskip
\cl{$\,\tilde X_S\,:=\;h\big(
\bigcup\,\big\{\,\bigcup\,\{\,[a,a+\epsilon(a)]\;|
\;a\in K_n\,\}\;\,\big|\,\;n\in S\,\big\}\big)\,\cup\,[1,2]\,$}
\smallskip 
for arbitrary infinite $\,S\subset{\Bbb N}\,$ then all spaces are 
order-isomorphic and hence homeomorphic!
(In fact, with the notation as in the proof of Proposition 2 below,
for any infinite $\,S\subset{\Bbb N}\,$ 
the linearly ordered set   $\,({\cal C}(\tilde X_S),\prec)\,$ 
is well-ordered and order-isomorphic 
to the set of all ordinal numbers $\,\alpha\leq \omega^\omega\,$.)

\bigskip
\medskip
{\bff 4. Proof of Theorem 3} 
\mp
The proof of Theorem 1 can be adapted in order to verify Theorem 3.
We replace each point set
$\;X_S\,=\,\bigcup\,\{\,[a_j,b_j]\;|\;j\in{\Bbb N}\,\}\;$ 
with $\;\tilde X_S\,=\,\bigcup\,\{\,[a_j,b_j]^n\;|\;j\in{\Bbb N}\,\}\;$
and claim that these compact subspaces of $\,{\Bbb R}^n\,$
are mutually non-homeomorphic also for arbitrary dimensions $\,n\,$.
As a consequence, Theorem 3 is settled by defining
$\,{\cal V}_n\,$ as the family of all open sets 
$\;\bigcup_{j=1}^\infty ]a_j,b_j[^n\;$
corresponding to $\,X_S\,$ represented as above with $\,S\,$
running through the infinite subsets of $\,{\Bbb N}\,$.
Indeed, (i), (ii), (iii) are obviously satisfied and (iv)  
follows from (v) since $\,\overline{V}\,$ is 
the completion of each metric space $\,V\in {\cal V}_n\,$.
\mp
Since an elimination of one point of a cube never destroys its 
connectedness,  
we cannot adopt the argumentation 
using the set $\,B_S\,$ in higher dimensions.
But fortunately we can stay very close to the proof of dimension 1
by transforming the concept of Cantor derivatives
from point sets of a topological space to
{\it families of components} of the space in the following way.
\mp
Let $\,{\cal G}\,$ be the family of all components of a Hausdorff space.
For every $\,{\cal F}\subset{\cal G}\,$ define 
$\;{\cal F}'\,=\,{\cal F}^{(1)}\,:=\;\{\,G\in{\cal G}\;|\;
G\cap\overline{\bigcup({\cal F}\setminus\{G\})}
\not=\emptyset\,\}\;$ and $\;{\cal F}^{(k+1)}:=({\cal F}^{(k)})'\;$
for every $\,k\in{\Bbb N}\,$. Now referring to $\,\tilde X_S\,$
let $\,\tilde{\cal L}_S\,$ be the family of all components 
$\,C\in{\cal G}\,$ such that $\,C\cap\overline{\tilde X_S\setminus C}\,$
contains precisely two points. Then, similarly 
as in the proof of Theorem 1, 
the set $\,S\,$ is topologically characterized via
\mp
\cl{$\;S\,=\,\big\{\,\min\,\{\,m\in{\Bbb N}\;|\;C\not\in {\cal G}^{(m+1)}\,\}
\;\big|\;C\in\tilde{\cal L}_S\,\big\}\;.$} 
\vfill\eject
\bigskip
\medskip
{\bff 5. Proof of Theorem 2} 
\mp
In the following, if $\,X\subset{\Bbb R}^n\,$ then 
$\,\overline X\,$ is the closure of $\,X\,$ in the space $\,{\Bbb R}^n\,$ 
and if $\,n=2\,$ then  $\,X^\circ\,$ is the interior of $\,X\,$ 
in the plane $\,{\Bbb R}^2\,$.
For abbreviation let $\,I=[0,1]\,$ and $\,J=\,]0,1[\;$ and let
$\;2{\Bbb N}\,:=\,\{\,2k\;|\;k\in{\Bbb N}\,\}\;$ be the set of all positive even numbers.
Furthermore, if $\;X\subset {\Bbb R}^2\;$ and $\;L\subset{\Bbb R}\;$
we regard $\,X\times L^k\,$ as a subset of $\,{\Bbb R}^{k+2}\,$
for every $\,k\geq 0\,$ where $\,X\times L^k\,$ is identified 
with $\,X\,$ if $\,k=0\,$.
\medskip
For each $\,m\in{\Bbb N}\,$ let $\;D_m\,:=\,[2^{-m},2^{-m+1}]^2\;$ and
\smallskip
\centerline{$W_m\;\,:=\,\;\bigcup\limits_{k=1}^m
\big[\,2^{-m}+{2k-1\over 2m+1}\cdot 2^{-m}\,,\,
2^{-m}+{2k\over 2m+1}\cdot 2^{-m}\,\big]^2\;.$}
\smallskip
So $\,D_m\,$ is a compact
square area and $\,W_m\,$ is a union of $m$ disjoint compact square areas
which all lie in the interior of $\,D_m\,$.
For $\;S\subset 2{\Bbb N}\;$ put
\smallskip
\centerline{$Z_S\;\,:=\,\;
[-1,0]^2\,\cup\bigcup\limits_{m\in S} (D_m\setminus W_m^\circ)\;$}
\smallskip
and define 
$\;{\cal F}_n\,:=\,
\{\,Z_S^\circ\times J^{n-2}\;|\;S\subset 2{\Bbb N}\,\}\;$
for each dimension $\,n\geq 2\,$.
Clearly, we always have $\,\overline{Z_S^\circ}=Z_S\,$
and hence $\;\overline{Z_S^\circ\times J^{n-2}}=Z_S\times I^{n-2}\,$.
Obviously, for every $\;U\in{\cal F}_n\;$ the closure
$\,\overline U\,$ is compact and a union of 
cubes $\,[a_1,a_1+h]\times\cdots\times[a_n,a_n+h]\,$ with $\,h>0\,$ so that 
the interiors of distinct cubes are always disjoint.
(Notice that $\,(D_m\setminus W_m^\circ)\times I^{n-2}\,$
is a union of precisely 
$\,((2m\!+\!1)^2-m)({1\over l})^{n-2}\,$
such cubes with edge length $\;l={2^{-m}\over 2m+1}\,$.)
\medskip
In order to verify that the family $\,{\cal F}_n\,$ 
has the desired homeomorphism properties
it is enough to investigate the components 
of the space $\,Z_S^\circ\times J^{n-2}\,$ and $\,Z_S\times I^{n-2}\,$ 
respectively. Clearly the components are always path connected
spaces and so it is natural to determine their fundamental groups.
(Two spaces $\,X,Y\,$ cannot be homeomorphic if
the fundamental group of some path component of $\,X\,$
is not isomorphic to the fundamental group of any path component of $\,Y\,$.)
\sp

For each $\;S\subset 2{\Bbb N}\;$ 
the components of the space
$\,Z_S\times I^{n-2}\,$ resp.~$\,Z_S^\circ\times J^{n-2}\,$ 
are precisely 
$\;(D_m\setminus W_m^\circ)\times I^{n-2}\;$ 
resp.~$\;(D_m^\circ\setminus W_m)\times J^{n-2}\;$ 
with $\;m\in S\;$ and the one simply connected component 
$\,[-1,0]^2\times I^{n-2}\,$ resp.~$\,]\!-1,0[^2\times J^{n-2}\,$.

For each $\;m\in{\Bbb N}\;$ 
the fundamental group both of $\;D_m\setminus W_m^\circ\;$
and of $\;D_m^\circ\setminus W_m\;$
is free on $\,m\,$ generators. This is enough since
for $\,n\geq 3\,$ both $\,I^{n-2}\,$ and $\,J^{n-2}\,$ 
have trivial fundamental groups. (If $\,X,Y\,$ are path connected
spaces then the fundamental group of the product space 
$\,X\!\times\! Y\,$ is isomorphic to the direct product of the
fundamental groups of $X$ and $Y$.)
\medskip
Finally, dispensing with compactness, it is plain to modify 
the definition of $\,{\cal F}_n\,$ so that 
each member of $\,{\cal F}_n\,$ is the interior of
a union of cubes of the form 
$\,[k_1,k_1\!+\!1]\times\cdots\times[k_n,k_n\!+\!1]\,$ with
$\;k_1,...,k_n\in{\Bbb Z}\,$. 
For example, for $\,\emptyset\not= S\subset 2{\Bbb N}\,$
replace $\,Z_S^\circ\times J^{n-2}\,$ 
with $\,Y_S^\circ\times J^{n-2}\,$ where 
$\;\,Y_S\;:=\;
\bigcup\limits_{m\in S} 
\{\,(t_m x,t_m y)\;|\;
(x,y)\,\in\,D_m\setminus W_m^\circ\,\}\;\,$
with $\,t_m:=4^m(2m+1)\,$.
\bigskip\mp
{\bff 6. Proofs of Propositions 1 and 3} 
\mp
First we need two basic lemmas. A proof of Lemma 1 is an easy exercise
and Lemma 2 is a consequence of a well-known theorem due to Sierpinski
(cf.~[1] 6.1.27).
\medskip
{\bf Lemma 1.} {\it If $\,g\,$ is an arbitrary injection from 
$\,{\Bbb N}\,$ into $\,{\Bbb N}\,$ then 
$\;\{\,n\in{\Bbb N}\;|\;n\leq g(n)\,\}\;$ must be 
an infinite set.}
\medskip
{\bf Lemma 2.} {\it If $\,a,b\in{\Bbb R}\,$ and $\,a<b\,$ then 
for any family $\,{\cal F}\,$ of mutually disjoint intervals 
$\,[x,y]\,$ with $\,x<y\,$ the equality 
$\;[a,b]\,=\,\bigcup{\cal F}\;$ is only possible 
in the trivial case $\,{\cal F}=\{[a,b]\}\,$.} 
\medskip\smallskip
In order to prove Proposition 1 it is enough to settle the statement 
on the closed point sets $\,X_u\,$ because each $\,X_u\,$ is
a completion of the metric space $\,X_u^\circ\,$.
Fix $\;2\leq v<w\;$
and let $\;I_n\,:=\,\big[2^{v^n},2^{v^n}+v^n\big]\,=:\,[a_n,b_n]\;$
and $\;J_n\,:=\,\big[2^{w^n},2^{w^n}+w^n\big]\,=:\,[c_n,d_n]\;$
for every $\,n\in{\Bbb N}\,$. 
So we have $\;b_n-a_n\,=\,v^n\;$
and  $\;d_n-c_n\,=\,w^n\;$ for every $\;n\in{\Bbb N}\;$
and $\;1+b_n\leq a_{n+1}\;$ and $\;1+d_n\leq c_{n+1}\;$
for every $\;n\in{\Bbb N}\,$.
Assume indirectly that $\,f\,$ is a uniformly continuous
bijection from $\,X_v\,$ onto $\,X_w\,$.
We claim that for each $\,n\in{\Bbb N}\,$ we must have 
$\,f(I_n)=J_{m}\,$ for some $\,m\in{\Bbb N}\,$.
Indeed, choose $\,m\,$ so that $\,f(I_n)\cap J_{m}\not=\emptyset\,$
and define an equivalence relation on $\,J_m\,$ via $\,x\sim y\,$ 
if and only if $\;f^{-1}(x),f^{-1}(y)\in I_k\;$ for some $\,k\,$. 
Then, since $\,f(I_k)\,$ is always compact and connected
and since all point sets $\,J_k\,$ are open and closed in the space $\,X_w\,$,
the family $\,{\cal F}\,$ of all equivalence classes 
must equal $\;\{\,f(I_k)\;|\;k\in K\,\}\;$ for some $\,K\subset{\Bbb N}\,$
with $\,n\in K\,$. Hence, in view of Lemma 2 we must have $\,K=\{n\}\,$ 
or, equivalently, $\,f(I_n)=J_{m}\,$.
\smallskip
Consequently, there is a
bijection $\;g\,:\;{\Bbb N}\to{\Bbb N}\;$ such that
$\;f(I_n)=J_{g(n)}\;$ for every $\,n\in{\Bbb N}\,$.
By Lemma 1, $\;G\,:=\,
\{\,n\in{\Bbb N}\;|\;n\leq g(n)\,\}\;$ is an infinite set.
Let $\;Y_v\,:=\,\bigcup\limits_{n\in G}I_n\;$ and 
$\;Y_w\,:=\,\bigcup\limits_{n\in G}J_{g(n)}\,$.
Then $\,f\,$ is a uniformly continuous function 
from the unbounded set $\,Y_v\,$ onto the unbounded set $\,Y_w\,$.
Thus we may fix $\,0<\delta<1\,$ so that 
$\;|f(x)-f(y)|\leq 1\;$ whenever $\;x,y\in Y_v\;$ and 
$\;|x-y|\leq \delta\,$.
Naturally, $\,f\,$ is strictly monotonic 
on each interval $\;I_n\,(n\in G)\;$ and 
$\;\{f(a_n),f(b_n)\}=\{c_{g(n)},d_{g(n)}\}\;$ for every $\,n\in G\,$.  
Now, for every $\,n\in G\,$ we have 

\centerline{$w^n\,\leq\,w^{g(n)}\,=
\,d_{g(n)}-c_{g(n)}\,=\,|f(a_n)-f(b_n)|\,=\,$}
\smallskip
\centerline{$|f(a_n)-f(a_n+\delta)|+|f(a_n+\delta)-f(a_n+2\delta)|+\cdots+
|f(a_n+k\delta)-f(b_n)|\,\leq\,k+1$}
\smallskip
where $\,k\in{\Bbb N}\,$ is chosen so that 
$\;a_n+k\delta<b_n\leq a_n+(k+1)\delta\;$ or, equivalently,
$\;k\delta <v^n\leq (k+1)\delta\,$. But then $\;w^n-1\leq v^n/\delta\;$
for every $\,n\,$ in the infinite set $\,G\,$.  
This is impossible since $\;\lim\limits_{n\to\infty}w^n/v^n\,=\,\infty\;$
and so the proof of Proposition 1 is finished.
\medskip\smallskip
{\it Remark.} 
Concerning higher dimensions, in view of the preceding proof it is plain
that the $\,c\,$ closed resp.~open point sets 
$\,\bigcup_{m=1}^\infty[2^{u^m},2^{u^m}+u^m]^n\,$ 
resp.~$\,\bigcup_{m=1}^\infty]2^{u^m},2^{u^m}+u^m[^n\,$ 
$\,(\,u\geq 2)\,$
in $\,{\Bbb R}^n\,$ are metrically distinct and topologically similar
for arbitrary $\,n\,$.
\bigskip
Now we are going to prove Proposition 3. As usual, 
the distance $\,d(A,B)\,$ between two nonempty subsets $\,A,B\,$ of 
$\,{\Bbb R}^n\,$ is the infimum of all numbers $\;d(a,b)\;$ with
arbitrary $\,a\in A\,$ and $\,b\in B\,$
where $\;d(x,y)\;$ denotes the Euclidean distance between 
$\,x,y\in{\Bbb R}^n$. 
If $\,U\,$ is an open subspace of $\,{\Bbb R}^n\,$
and if $\,{\cal G}\,$ is the (countable) family of 
all components of $\,U\,$, then 
let us call a finite subset $\,{\cal F}\,$ of $\,{\cal G}\,$
a {\it chain} if and only if there is an ordering 
$\;{\cal F}\,=\,\{U_1,...,U_m\}\;$ with $\,U_i\not=U_j\,(i\not=j)\,$
and $\;d(U_k,U_{k+1})=0\;$ for every $\,k<m\,$.
The {\it length} of $\,{\cal F}\,$ is $\,m\,$.
A chain is {\it maximal} if it is not contained in a chain of greater length.
For every $\,S\subset{\Bbb N}\,$ all the components of 
$\,Y_S\,$ are open cubes of the form 
$\;]a_1,a_1+h[\,\times\cdots\times\,]a_n,a_n+h[\;$ 
and, evidently (by induction on the dimension $n$),
\smallskip
\centerline{$Y_S\;=\;
\bigcup\limits_{s\in S} \big(\bigcup{\cal F}_{s}\big)\;\cup\; 
\bigcup\limits_{k=1}^\alpha C_k\qquad
(\,\alpha\,\in\,{\Bbb N}\cup\{\infty\}\,)$}
\smallskip
where $\,\{C_k\}\,$ is always a maximal chain of length $\,1\,$
and $\,{\cal F}_{s}\,$ is a maximal chain of length $\,(s+1)^n\,$
for every $\,s\in S\,$ and 
$\;d(\bigcup{\cal F}_i,\bigcup{\cal F}_j)>0\;$ 
whenever $\,i,j\in S\,$ and $\,i\not=j\,$.
\medskip
Suppose that $\,S,S'\subset{\Bbb N}\,$ 
and that $\,f\,$ is a uniform homeomorphism from 
$\,Y_S\,$ onto $\,Y_{S'}\,$. 
Of course, $\,f(C)\,$ is a component of $\,Y_{S'}\,$   
if and only if $\,C\,$ is a component of $\,Y_{S}\,$.   
Furthermore, for any
$\,\emptyset\not=A,B\subset Y_{S}\,$ we certainly have
$\,d(A,B)=0\,$ if and only if $\,d(f(A),f(B))=0\,$.
Therefore, for every $\,s\in S\,$ the set 
$\;\{\,f(C)\;\,|\,\;C\in {\cal F}_s\,\}\;$ 
must be a maximal chain of length $\,(s+1)^n\,$ 
in the space $\,Y_{S'}\,$, whence $\,s\in S'\,$.
Thus $\,S\subset S'\,$. Similarly, $\,S'\subset S\,$.
\bp\mp
{\bff 7. Proof of Proposition 2}
\mp
For a nonempty set $\,X\,$ 
in the family $\,{\cal U}\cup{\cal A}\,$
let $\,{\cal C}(X)\,$ be the family of all 
components of the Euclidean subspace $\,X\,$ of $\,{\Bbb R}\,$.
Since each member of $\,{\cal C}(X)\,$
is an open or closed interval, we may define
a natural strict linear ordering $\,\prec\,$ of $\,{\cal C}(X)\,$
via $\,A\prec B\,$ for distinct (and hence disjoint)
$\,A,B\in{\cal C}(Z)\,$
if and only if $\,a<b\,$ for {\it some} $\,(a,b)\in A\times B\,$
or, equivalently,  
if $\,a<b\,$ for {\it every} $\,(a,b)\in A\times B\,$.
\mp
Let $\;\emptyset\not=X,Y\subset{\Bbb R}\;$ 
and let $\;\varphi:\,X\to Y\;$ be an order isomorphism. 
Then $\,\varphi\,$ is a homeomorphism with respect 
to the order topologies of $\,(X,<)\,$ and $\,(Y,<)\,$.
Moreover, if the sets $\,X,Y\,$ lie in the family
$\,{\cal U}\cup{\cal A}\,$ then 
$\,\varphi\,$ is a homeomorphism between the Euclidean spaces
$\,X\,$ and $\,Y\,$ and hence 
$\;A\mapsto \varphi(A)\;$ defines a bijection 
from $\,{\cal C}(X)\,$ onto $\,{\cal C}(Y)\,$
and it is evident that this bijection is an order isomorphism 
between $\,({\cal C}(X),\prec)\,$ and $\,({\cal C}(Y),\prec)\,$.
Thus, $\;X,Y\,\in\,{\cal U}\cup{\cal A}\;$ 
are not order-isomorphic if 
the two families $\,{\cal C}(X)\,$ and $\,{\cal C}(Y)\,$
are not order-isomorphic.
(Conversely, if either $\,X,Y\in {\cal U}\,$ or $\,X,Y\in {\cal A}\,$ 
are compact then
from any order isomorphism between 
$\,({\cal C}(X),\prec)\,$ and $\,({\cal C}(Y),\prec)\,$
we may easily construct an order isomorphism between 
$\,(X,<)\,$ and $\,(Y,<)\,$.
This is not true for arbitrary sets $\,X,Y\in {\cal A}\,$
or for compact sets $\,X,Y\subset{\Bbb R}\,$.
Consider, for example, 
$\;X\,=\,[0,1]\cup[2,3]\;$ and firstly $\;Y\,=\,[0,1]\cup[2,\infty[\,$
and secondly $\;Y\,=\,[0,1]\cup\{2\}\,$.) 
\mp
So in order to settle Proposition 2 it is enough to find 
$\,c\,$ metrically similar sets $\,X\,$
in the family $\,{\cal U}\,$ resp.~$\,{\cal A}\,$,
such that the corresponding sets $\,{\cal C}(X)\,$
are mutually not order-isomorphic. 
Let $\,{\cal G}\,$ be the family of 
all functions $\,g\,$ from $\,{\Bbb N}\,$ to $\,\{0,1\}\,$ 
such that the set $\,g^{-1}(\{1\})\,$ is infinite.
Clearly, the cardinal number of the family $\,{\cal G}\,$ is $\,c\,$.
For every $\,n\in{\Bbb N}\,$ define 
\sp
\cl{$\,Z_n\;:=\;\bigcup\limits_{k=1}^\infty
\big(\,]6n+1+3^{-2k},6n+1+3^{-2k+1}[\;\cup\; 
]6n+2-3^{-2k+1},6n+2-3^{-2k}[\,\big)$}
\sp
and for every $\,g\in{\cal G}\,$ define
\sp
\cl{$\,U_g\;:=\;\bigcup\limits_{n=1}^\infty 
\big(\,Z_n\cup\,]6n,6n+1[\,\cup\,]6n+2,6n+3[\,\big)
\;\cup\;\bigcup\limits_{n\in g^{-1}(\{1\})}]6n+4,6n+5[\;.$}
\sp
It is plain that the $\,c\,$ open sets $\;U_g\,(g\in{\cal G})\;$ 
are metrically similar and that  
the $\,c\,$ closed sets $\;\overline{U_g}\,(g\in{\cal G})\;$ 
lie in the family $\,{\cal A}\,$ and are  metrically similar too.
Let $\,\zeta\,$ denote the order type of $\,{\Bbb Z}\,$.
Then $\,\zeta\,$ is also the order type of $\,{\cal C}(Z_n)\,$
for every $\,n\in{\Bbb N}\,$
and it is evident that 
for each $\,g\in{\cal G}\,$ both the order type of $\;{\cal C}(U_g)\;$ 
and the order type of $\;{\cal C}(\overline{U_g})\;$
equals
\mp
\cl{$(1+\zeta+1)\,+\,g(1)\,+\,
(1+\zeta+1)\,+\,g(2)\,+\,
(1+\zeta+1)\,+\,g(3)\,+\,\cdots\,$\quad}
\sp
\cl{$=\;1\,+\,\zeta\,+\,(1+g(1)+1)\,+\,\zeta\,+\,(1+g(2)+1)
\,+\,\zeta\,+\,(1+g(3)+1)\,+\,\cdots\,$}
\mp
where a nonnegative integer $\,k\,$
is always the order type of any linearly ordered set 
of precisely $\,k\,$ elements.
(If $\,\alpha,\beta\,$ are order types of nonempty sets
then $\,\alpha+0+\beta\,$ is just $\,\alpha+\beta\,$.)
Naturally, for distinct $\,f,g\in{\cal G}\,$ the order types 
of $\,{\cal C}(U_f)\,$ and $\,{\cal C}(U_{g})\,$
are distinct and this concludes the proof.
\mp
{\it Remark.} In view of the previous considerations 
it is easy to track down 
$\,c\,$ metrically similar compact subsets of $\,{\Bbb R}\,$
without singleton components which are mutually not order-isomorphic. 
(The existence of such sets follows from [5] Main Theorem 2.)
Take for example  (with $\,h\,$ as in the proof of Theorem 1)
the $\,c\,$ point sets $\;h(\overline{U_g})\cup [1,2]\;(g\in{\cal G})\,$.  
\vfill\eject
\bp\sp
{\bff 8. Totally disconnected point sets}
\mp
So far we considered only point sets where no component 
is a singleton. Now we consider point sets 
where every component 
is a singleton, i.e.~{\it totally disconnected} point sets.
Since for every totally disconnected set $\,T\subset{\Bbb R}\,$
the set $\,T\times\{0\}^{n-1}\,$ is a
totally disconnected subset of the Euclidean space $\,{\Bbb R}^n\,$,
for our purpose
there will be no benefit of considering arbitrary dimensions 
and so we restrict to dimension $\,n=1\,$ in the following.
(In view of Theorem 4 below, notice also that every compact and
totally disconnected subspace of $\,{\Bbb R}^n\,$
is homeomorphic to some subspace of $\,{\Bbb R}\,$.)
In the real line $\,{\Bbb R}\,$ a point set
is totally disconnected 
if and only if it does not contain a nondegenerate interval.
(In particular, no nonempty open set is totally disconnected.)
The real line $\,{\Bbb R}\,$ contains $\,2^c\,$ totally disconnected
subsets, for example all subsets of $\,{\Bbb R}\setminus{\Bbb Q}\,$.
Among these $\,2^c\,$ totally disconnected spaces there must also be
$\,2^c\,$ non-homeomorphic spaces because one cannot track down 
more than $\,c\,$ homeomorphic subspaces
of $\,{\Bbb R}\,$.
(For if 
$\,X\subset{\Bbb R}\,$ then there are at most $\,c\,$ continuous functions
from $\,X\,$ into $\,{\Bbb R}\,$.) 
\sp
What is the number of all topological types of 
totally disconnected {\it closed} point sets in $\,{\Bbb R}\,$?
The following theorem, 
which is a noteworthy counterpart to 
Theorem 1, gives the answer.
\mp
{\bf Theorem 4.} {\it There are $\,c\,$ mutually non-homeomorphic 
compact and totally disconnected 
subspaces of the Euclidean unit interval $\,[0,1]\,$.}
\mp
{\it Proof.} In any Hausdorff space $\,X\,$
a point set $\,A\,$ is {\it dense in itself}
if and only if every point in $\,A\,$ is a limit point of $\,A\,$,
i.e.~$\,A\subset A'\,$.
Let $\;\Delta(X)\,:=\,\bigcup\,\{\,A\subset X\;|\;A\subset A'\,\}\;$
denote {\it the maximal dense-in-itself point set} in the space $\,X\,$.
Define a sort of {\it signature} set of integers by
\sp
\centerline{$\;\;\Sigma(X)\;:=\;
\big\{\,k\in{\Bbb N}\;\,\big|\,\;
\big((X\setminus \Delta(X))^{(k)}\setminus(
X\setminus \Delta(X))^{(k+1)}\big)\cap \Delta(X)  
\,\not=\,\emptyset\,\big\}\;.$}
\sp
Let $\,h\,$ be a strictly increasing function from $\,[0,\infty[\,$
onto $\,[0,1[\,$, for example $\,h(x)={2\over\pi}\arctan x\,$
as in the proof of Theorem 1.
Let $\,{\Bbb D}\subset [0,1]\,$ be the Cantor ternary set.
(Notice that $\,{{\Bbb D}}\;\!'={\Bbb D}\,$.)
As in the proof of Theorem 1,
for every $\,n\in{\Bbb N}\,$ let 
$\,K_n\subset[5n,5n+1]\,$ be compact with $\;\max K_n\,=\,5n+1\;$ 
such that the $k$-th derivative $\,K_n^{(k)}\,$ is infinite 
whenever $\,k<n\,$ and empty whenever $\,k>n\,$ and a singleton 
when $\,k=n\,$, namely $\;K_n^{(n)}=\{5n+1\}\,$.
\sp
Now for each infinite set $\,S\subset{\Bbb N}\,$ define 
\sp
\cl{$D_S\;:=\;\bigcup\limits_{n\in S}\{\,5n+1+x\;\,|\,\;x\in{\Bbb D}\,\}\;$\qquad 
and\qquad$\;Y_S\;:=\;
h\,\big(\,D_S\,\cup\bigcup\limits_{n\in S} K_n\,\big)\,\cup\,\{1\}\;.$}
\sp
Of course, all point sets $\,Y_S\subset [0,1]\,$ 
are compact and totally disconnected.
Moreover, $\;\Delta(Y_S)\,=\,h(D_S)\cup\{1\}\;$ and
$\;(K_n\setminus D_S)^{(k)}=K_n^{(k)}\;$ for every $\,k\in{\Bbb N}\,$ and 
$\,n\in S\,$. Therefore, similarly as in the proof of [3] Theorem 7.1
we always have $\;\Sigma(Y_S)=S\;$ and hence 
$\,Y_S\,$ and $\,Y_{T}\,$ are never homeomorphic 
for distinct infinite sets $\,S,T\subset{\Bbb N}\,$, {\it q.e.d.}
\mp\sp
{\it Remark.} If in Theorem 4 the property {\it perfect} 
({\it dense-in-itself}) is 
added then the variety of the spaces collapses.
Indeed, it is well-known that
any perfect, compact, zero-dimensional, second countable 
Hausdorff space is homeomorphic to $\,{\Bbb D}\,$ (cf.~[2]).
(Note that a compact Hausdorff space is zero-dimensional 
if and only if it is totally disconnected, cf.~[1] 6.2.9). 
\mp\sp
{\it Remark.} By a classic theorem due to Mazurkiewicz und Sierpinski [4] 
there are precisely $\,\aleph_1\,$ compact and countable
Hausdorff spaces up to homeomorphism. 
($\,\aleph_1\,$ is the least cardinal number greater than
the cardinality $\,\aleph_0\,$ of a countably infinite set,
whence $\,\aleph_1\leq c\,$.)
As a consequence, 
since each countable metric space can be embedded in $\,{\Bbb R}\,$
(see [1] 4.3.H.b),
the space $\,{\Bbb R}\,$
has uncountably many 
non-homeomorphic compact and totally disconnected 
subspaces.
Theorem 4 is an improvement of this consequence 
because one cannot rule out the existence 
of an uncountable set whose cardinality is smaller than $\,c\,$.
(Actually, it is consistent with ZFC set theory that 
there exist $\,c\,$ uncountable cardinal numbers smaller than $\,c\,$.)
\bp
{\bff 9. Discrete linear point sets}
\mp
There is an interesting consequence of Theorem 4 concerning {\it discrete}
point sets.
A point set $\,X\,$ in the real line 
is discrete if and only if for every $\,x\in X\,$
the singleton $\,\{x\}\,$ is open in the 
subspace $\,X\,$ of $\,{\Bbb R}\,$.
Equivalently, for every $\,x\in X\,$ there is 
$\,\delta_x>0\,$ such that 
$\,[y-\delta_y,y+\delta_y]\cap [z-\delta_z,z+\delta_z]=\emptyset\,$ 
whenever $\,y,z\in X\,$ and $\,x\not=y\,$.
Consequently, any discrete subset of $\,{\Bbb R}\,$ is countable.
Since $\,{\Bbb R}\,$ has precisely $\,c\,$ countable subsets 
and since $\,x+{\Bbb Z}\,$ is discrete for $\,0<x<1\,$,
the real line contains precisely $\,c\,$
discrete point sets. 
From the topological point of view, essentially there is 
precisely one infinite discrete point set in $\,{\Bbb R}\,$.
Indeed, if $\,X\,$ is any infinite discrete point  
set in $\,{\Bbb R}\,$ then $\,X\,$ is obviously homeomorphic to 
the discrete space $\,{\Bbb Z}\,$. 
On the other hand,  from the metrical point of view 
there are very many discrete point sets.
\mp
{\bf Theorem 5.} {\it There are $\,c\,$ metrically distinct
infinite and discrete point sets in the unit 
interval $\,[0,1]\,$.} 
\mp
{\it Proof.} Choose for every set $\,Y_S\,$
in the proof of Theorem 4
a discrete set $\,Z_S\subset [0,1]\,$ 
such that $\;Z_S\cap Y_S\,=\,\emptyset\;$
and $\;Z_S'=Y_S\,$, whence $\,Z_S\,$ is infinite and
$\;\overline{Z_S}\setminus Z_S\,=\,Y_S\,$.
Suppose that for infinite sets $\,S,T\subset{\Bbb N}\,$ there is 
a bijection $\,f\,$ from $\,Z_S\,$ onto $\,Z_T\,$ such that
$\,f,f^{-1}\,$ are uniformly continuous.
Then there is an expansion of $\,f\,$ to a homeomorphism $\,g\,$
from $\,\overline{Z_S}\,$ onto $\,\overline{Z_T}\,$ 
since the Euclidean metric space 
$\,\overline{Z_S}\,$ resp.~$\,\overline{Z_T}\,$   
is a completion of the Euclidean metric space 
$\,Z_S\,$ resp.~$\,Z_T\,$. Then we must have $\,S=T\,$ since 
$\;g(Y_S)\,=\,g(\overline{Z_S}\setminus Z_S)\,=\,
g(\overline{Z_S})\setminus f(Z_S)\,=\,
\overline{Z_T}\setminus Z_T\,=\,Y_T\;$
and $\,Y_S,Y_T\,$ cannot be homeomorphic if $\,S\not=T\,$.
It is always possible to choose such sets $\,Z_S\,$ 
and the proof of Theorem 5 is finished in 
view of the following lemma.
\mp
{\bf Lemma 3.} {\it Any compact and totally disconnected and nonempty set
$\,A\subset{\Bbb R}\,$ equals $\,Z'\,$ for some discrete set 
$\,Z\subset[\min A,\max A]\,$
with $\,Z\cap A=\emptyset\,$.}
\mp
{\it Proof.} Let $\,a=\min A\,$ and $\,b=\max A\,$.
Clearly $\,A\,$ must be nowhere dense. Hence 
$\,A\,$ is the boundary of the open set $\,[a,b]\setminus A\,$.
Write $\;[a,b]\setminus A\,=\,\bigcup_{j\in J}I_j\;$
with a countable index set $\,J\,$ and where the sets $\,I_j\,$ are mutually
disjoint open intervals, 
$\;I_j=\,]x_j,y_j[\;$ with $\,x_j<y_j\,$. 
Define  a countable 
set $\,C_j\subset I_j\,$ by
\sp
\cl{$\;C_j\,:=\,\{\,x_j+(y_j-x_j)2^{-k}\;|\;k\in{\Bbb N}\}\cup
\{\,y_j-(y_j-x_j)2^{-k}\;|\;k\in{\Bbb N}\}\;$}
\sp
for every $\,j\in J\,$
and put $\;Z\,=\,\bigcup_{j\in J}C_j\,$.
Clearly, $\,Z\,$ is discrete and $\,Z\subset [a,b]\,$ 
and $\,Z\cap A=\emptyset\,$.
Since  $\,C_j'=\{x_j,y_j\}\,$ for every $\,j\in J\,$,
we have $\,Z'=A\,$, {\it q.e.d.}
\bp
Since the discrete sets in Theorem 5 are all {\it absolutely bounded},
none of them is closed. 
Naturally, a closed and discrete subset of $\,{\Bbb R}\,$ cannot be bounded
and infinite, while
a bounded and discrete subset of $\,{\Bbb R}\,$ cannot be infinite and closed.
Obviously, an infinite closed set $\,A\subset{\Bbb R}\,$ is discrete if and only if 
$\,A\cap [-k,k]\,$ is finite for every $\,k\in{\Bbb N}\,$.
How many infinite, closed, discrete point sets do exist
from the metrical point of view?
The following theorem gives the answer.
\mp
{\bf Theorem 6.} {\it There are $\,c\,$ metrically distinct
infinite, closed, discrete point sets in 
the real line $\,{\Bbb R}\,$.}
\mp
{\it Proof.} Let $\,{\Bbb P}\,$ be the set of all primes.
For $\,p\in{\Bbb P}\,$ and $\,n\in{\Bbb N}\,$ define a set of precisely 
$\,p\,$ elements which is contained
in an interval of length $\,p^{-n}\,$ by 
\sp
\cl{$\;G[p;n]\;:=\;\{\,p^n+k^{-1}p^{-n}\;\,|\,\;k=1,..,p\,\}\,.$}
\sp
\eject
For $\,\emptyset\not=S\subset{\Bbb P}\,$ define $\,A_S\subset{\Bbb R}\,$ 
by $\;\;A_S\,:=\bigcup\limits_{p\in S}\bigcup\limits_{n=1}^\infty G[p;n]\,.$ 

We always have $\,G[p;n]\subset\,]p^n,p^n+1[\;$ and,
clearly, $\;]p^n,p^n+1[\,\cap\,]q^m,q^m+1[\,\,=\,\emptyset\;$ 
whenever $\,p,q\in{\Bbb P}\,$
and $\,n,m\in{\Bbb N}\,$ and $\,(p,n)\not=(q,m)\,$.
Thus $\,A_S\,$ is always a
closed and discrete point set. 
We claim that there is no uniform homeomorphism 
between any two of the $\,c\,$ Euclidean metric spaces 
$\;A_S\;(\emptyset\not=S\subset{\Bbb N})\,$. 
Let $\,S,T\subset{\Bbb P}\,$ and $\,S\not\subset T\,$
and suppose indirectly that $\,f\,$ is a bijection from $\,A_S\,$ onto
$\,A_T\,$ such that $\,f\,$ and $\,f^{-1}\,$ are uniformly continuous.
Then we can choose $\,\delta>0\,$
so that 
\sp
(i) $\;\;|f(x)-f(y)|<{1\over 2}\;$ whenever $\,x,y\in A_S\,$
and $\,|x-y|<\delta\,$, 
\sp
(ii)
$\;|f^{-1}(x)-f^{-1}(y)|<{1\over 2}\;$ whenever $\,x,y\in A_T\,$
and $\,|x-y|<\delta\,$. 
\sp
Now choose $\,p\,$ in $\,S\setminus T\,$
and fix $\,N\in{\Bbb N}\,$ so that $\,p^{-N}<\delta\,$.
In view of (i), for every $\,n\geq N\,$ we can define 
a prime $\,q_n\in T\,$ and a number $\,m_n\in{\Bbb N}\,$
such that $\;f(G[p;n])\subset G[q_n;m_n]\,$.
Since 
$\,G[q_n;m_n]\,$ is finite,
not both $\;\{\,q_n\;|\;n\geq N\,\}\;$
and $\;\{\,m_n\;|\;n\geq N\,\}\;$ can be finite sets.

Hence we may choose $\,n\geq N\,$ so that $\,q_n^{-m_n}<\delta\,$.
But then, in view of (ii) and $\,f(G[p;n])\subset G[q_n;m_n]\,$
we have $\,f^{-1}(G[q_n;m_n])=G[p;n]\,$ and
this is impossible since 
$\,G[q_n;m_n]\,$ has precisely $\,q_n\,$ elements 
and $\,G[p;n]\,$ has precisely $\,p\,$ elements 
and $\,q_n\not=p\,$, {\it q.e.d.}
\bp
{\it Remark.} From the order-theoretical point of view
there are only countably many 
closed and discrete subsets of $\,{\Bbb R}\,$.
Firstly, two finite sets are order-isomorphic 
if and only if they are equipollent. 
Secondly, a moment's reflection is sufficient to see that
any infinite, closed, discrete subset of $\,{\Bbb R}\,$
is order-isomorphic 
to $\,{\Bbb Z}\,$ or to $\,{\Bbb N}\,$ or to $\,{\Bbb Z}\setminus{\Bbb N}\,$.
On the other hand, there are $\,c\,$ mutually not
order-isomorphic discrete subsets of $\,{\Bbb R}\,$.
For example, let $\;U_g, {\cal C}(U_g)\;(g\in{\cal G})\;$ be as in the proof 
of Proposition 2 and let $\,\varphi\,$ denote any choice 
function on the family $\;{\cal U}\setminus\{\emptyset\}\,$,
i.e.~$\;\varphi(U)\in U\;$ whenever $\;\emptyset\not= U\in {\cal U}\,$.
Then for each $\,g\in{\cal G}\,$ the set 
$\;D_g\,:=\,\{\,\varphi(C)\;\,|\,\;C\in{\cal C}(U_g)\,\}\;$
is discrete and the order type of $\,D_g\,$ equals 
the order type of $\,{\cal C}(U_g)\,$, whence 
$\,D_f\,$ and $\,D_g\,$ are never order-isomorphic for 
distinct $\,f,g\in{\cal G}\,$.
\mp
\bp
In a natural way the proof of Theorem 5 leads to
the following noteworthy theorem. 
\mp
{\bf Theorem 7.} {\it A countably infinite discrete space $\,X\,$ has $\,c\,$ 
mutually non-homeomorphic 
metrizable compactifications of size $\,c\,$.
There exist $\,c\,$ incomplete metrics $\,d\,$ 
on a countably infinite set $\,X\,$
such that $\,(X,d)\,$ is always a discrete topological space
and the completions of the metric spaces $\,(X,d)\,$ are compact of size 
$\,c\,$
and topologically distinct.}  
\mp
{\it Proof.} For each infinite set $\,S\subset{\Bbb N}\,$ let $\,Y_S,Z_S\,$ 
be as in the proof
of Theorem 5 and 
define a metric $\,d_S\,$ on $\,X\,$ such that the metric space $\,(X,d_S)\,$
is an isometric copy of $\,Z_S\,$ equipped with the Euclidean metric. 
Then the topology of the discrete space $\,X\,$ is induced by
the metric $\,d_S\,$ which of course is not complete.
The Euclidean metric space $\;\overline{Z_S}\,=\,Y_S\cup Z_S\;$
is both a completion of the metric space $\,(X,d_S)\,$
and a compactification of the discrete space $\,X\,$.
Two spaces $\;\overline{Z_S},\overline{Z_T}\;$ are never
homeomorphic for distinct infinite sets $\,S,T\subset{\Bbb N}\,$
since in view of $\,Z_S'=Y_S\,$ we have 
$\;\Sigma(\overline{Z_S})=\Sigma(Y_S)=S\;$ for each $\,S\,$, {\it q.e.d.}
\bp
{\it Remark.} The cardinality $\,c\,$ in Theorem 7 is  
the largest possible in all cases. Indeed, 
there exist precisely $\,c\,$ compact metrizable spaces  up to homeomorphism
and any compact metric space is either countable or of size $\,c\,$
(cf.[3]). And in view of [4] a countably infinite discrete space
has only $\,\aleph_1\,$ topologically distinct 
countable (and hence metrizable) compactifications. 
(Note that if $\,\alpha\,$ is any countable ordinal number 
then in the compact space 
$\,\omega^\alpha+1\,$ the successor ordinals form a discrete 
and dense subspace.)
\bp\bp\bigskip
{\bff References}
\medskip
[1] Engelking, R.: {\it General Topology}, revised and completed edition.
Heldermann 1989. 
\smallskip
[2] Kechris, A.: {\it Classical Descriptive Set Theory}. Springer 1995.
\smallskip
[3] Kuba, G.: {\it Counting metric spaces}.
Arch.d.Math. {\bf 97}, 569-578 (2011).
\smallskip
[4] Mazurkiewicz, S.; Sierpinski, W.: Contribution … la topologie
des ensembles d'nombrables. 

\rightline{{\it Fund.~Math.}~{\bf 1} (1920), 17-27.}
\sp
[5] Winkler, R.: {\it How much must an order theorist forget to become a 
topologist?} 

\rl{Contributions to General Algebra 12, Proc.~Vienna Conference (June 1999).}

\bigskip\bp

{\sl Author's address:} Institute of Mathematics, 

University of Natural Resources and Life Sciences, Vienna, Austria. 
\smallskip
{\sl E-mail:} {\tt gerald.kuba@boku.ac.at}
\bp\bp\bp
\hrule\bp\bp\bp
{\bf Additional References}
\mp
[A1] Kuba, G.: {\it  On the variety of Euclidean point sets}.
Int.~Math.~News {\bf 228}, 23-32 (2015).  
\sp
[A2] Kuba, G.: {\it Counting ultrametric spaces}.
Colloq.~Math.~{\bf 152}, 217-234 (2018).
\bp\mp
{\sl Up to Section 7 the present paper is essentially identical with} [A1].
\mp
{\sl The proof of Theorem 5 is very similar to the proof of Proposition 2 in}
[A2]. {\sl Theorem 7 is a consequence of} [A2] Corollary 2.
{\sl In connection with Theorem 7 the following consequence of} 
[A2] Corollary 3 {\sl is worth mentioning.}
\mp
{\bf Theorem 8.} {\it The topology of an infinite discrete 
space $\,S\,$ can be generated by 
$\,2^{|S|}\,$ metrics $\,d\,$ such that 
the completions of the metric spaces $\,(S,d)\,$ are 
mutually non-homeomorphic metric spaces of size $\,|S|\,$.}
\mp
The cardinality $\,2^{|S|}\,$ in Theorem 8 is  
the largest possible since, trivially, an infinite 
set $\,S\,$ cannot carry more metrics than the total amount of 
mappings from $\,S\times S\,$ to $\,{\Bbb R}\,$ which equals 
$\,|{\Bbb R}|^{|S\times S|}=2^{|S|}\,$.
\end